\documentclass[11pt,twoside,leqno]{aomamlt2e}

\pdfoutput=1

\begin{document}
\title{{\center Integration of $e^{x^{n}}$ and $e^{-x^{n}}$ in forms of series, their applications in the field of differential equation; introducing generalized form of Skewness and Kurtosis; extension of starling's approximation.}}

\acknowledgements{T.S.K wants to thank Mr. Subhajit Paul, Department of Botany, University of Kalyani, Kalyani, India for his invaluable suggestions. S.P. wants to thank ``Deutsche Forschungsgemeinschaft'' (German Research Foundation) [DFG](GK 1177) for funding.}

\twoauthors{T. S. Konar $^1$}{   S. Paul $^2$}

           { \center $^1$ Jadavpur University, Kolkata, India  \\

            \center $^2$ Institute for Theoretical physics and Astrophysics, \\
            \center Wuerzburg University, Germany \\
            \hspace{4cm} \email{spaul@physik.uni-wuerzburg.de }}

\shorttitle{Integration of $ e^{x^n} $ and $ e^{-x^n} $ in forms of series}

\begin{abstract}
In this paper we tried a different approach to work out the integrals of $e^{x^{n}}$ and $e^{-x^{n}}$. Integration by parts shows a nice pattern which can be reduced to a form of series. We have shown both the indefinite and definite integrals of the functions mentioned along with some essential properties e.g. conditions of convergence of the series. Further more, we used the integrals in form of series to find out series solution of differential equations of the form $ x\frac{{d^2 y}}{{dx^2}}-(n-1)\frac{{dy}}{{dx}}-n^2 x^{2n-1}y-nx^n=0 $ and $ x\frac{{d^2 y}}{{dx^2}}-(n-1) \frac{{dy}}{{dx}}-n^2x^{2n-1}y+(n-1)=0 $, using some non standard method. We introduced modified Normal distribution incorporating some properties derived from the above integrals and defined a generalized version of Skewness and Kurtosis. Finally we extended Starling's approximation to $\mathop {\lim }\limits_{n \to \infty } (2n)! \sim 2n\sqrt {2\pi } \left( {\frac{{2n}} {e}} \right)^{2n}$
\end{abstract}

\section{Introduction}

Integrals of exponential functions like $e^{x^n}$ or $e^{-x^n}$ can generally be expressed in summation form $\sum\limits_{m = 0}^\infty= {\frac{{x^{nm +1}}}{{m!(nm+1)}}} + k_1 $ and  $ \sum \limits_{m = 0}^\infty {( - 1)^m \frac{{x^{nm + 1}}} {{m!(nm + 1)}}} + k_2 $ respectively. We tried a different approach to evaluate the integrals which can be applied to solve some special kind of mathematical problems. When we tried to integrate the functions mentioned, by the method of integration by parts, we see a nice pattern which can be reduced to a form of a series.

\newpage
In this paper we will show.
\begin{itemize}

\item For all positive  integer values of n, indefinite integral of $e^{x^n } $ and  $e^{- x^n }$ is \\
 $e^{x^n} \left[ {\sum\limits_{r = 0}^\infty  {\frac{{( - 1)^r n^r x^{1 + nr} }} {{\prod\limits_{p = 0}^r {(1 + pn)} }}} } \right] + C_1$ and $e^{ - x^n } \left[ {\sum\limits_{r = 0}^\infty  {\frac{{n^r x^{1 + nr} }} {{\prod\limits_{p = 0}^r {(1 + pn)} }}} } \right] + C_2$ respectively [ $\forall x;$ at least $\forall n \in \mathbb{Z}^ + $].\\

\item We showed the convergence of the above series to work out the proper integrals (e.g. 1.2.1) and shown that the modified interval of convergence for improper integrals are $ ]-\infty,\infty[ $ for n=even. And  $ -\infty \leqslant x < \infty $ for n=odd for $\int{e^{x^{n}}}dx $ and $[-\infty, \infty]$ for n=even and $-\infty<x\leqslant\infty $ for n=odd for $\int{e^{-x^n}}dx $.\\

\item We applied these to solve some special kind of differential equation of the form $ \;\;\;\; x\frac{{d^2y}}{{dx^2}}-(n-1)\frac{{dy}}{{dx}}-n^2x^{2n-1}y-nx^n=0 \;\;\;\; $   and   $\;\;\;\; x\frac{{d^2 y}}{{dx^2}}-(n-1)\frac{{dy}}{{dx}}-n^2 x^{2n-1}y+(n-1)=0 \;\;\; $ as series using some non standard method.

\item Further we used some properties of the above integrals to define some probability distribution function (P.D.F) and finally we defined a generalized version of Skewness, Kurtosis and Kurtosis excess for generalized "$(\frac{n}{2})^{th}$ order Normal distribution"  i.e. $y(x) = \frac{1}{{\sqrt[n]{n}\sigma P_n }}e^{ - \frac{1}{n}\left( {\frac{{x - m}}{\sigma }} \right)^n } $ for even n,\\

\begin{itemize}

\item Coefficient of Skewness = $ \sqrt[n]{{\frac{{m_{n + 1}^n }} {{m_n^{n + 1} }}}}$.

\item Kurtosis = $ \frac{{m_{2n} }}{{m_n^2 }}$.

\item Kurtosis excess = $\frac{{m_{2n} }} {{m_n^2 }} - (1 + n)$.

Where $m_s$  denotes $s_{th}$ order central moment of the distribution.

\end{itemize}

\item We also extended Starling's approximation to $\mathop {\lim }\limits_{n \to \infty } (2n)! \sim 2n\sqrt {2\pi } \left( {\frac{{2n}} {e}} \right)^{2n}$.

\end{itemize}

\section{Integration of $e^{x^{n}}$ and $e^{-x^{n}}$}

Generally integrals of $e^{x^{n}}$ and  $ e^{- x^{n}}$ can be
expressed in an easy summation form $ \sum\limits_{m = 0}^\infty
{\frac{{x^{nm +1}}}{{m!(nm+1)}}} + k_1 \;\;\;\; ...(1) $     and
     $ \sum\limits_{m = 0}^\infty  {( - 1)^m \frac{{x^{nm + 1}}}{{m!(nm + 1)}}}  + k_2 \;\;\;\; ...(2) $ respectively.But we tried to evaluate them in a different manner, which is rather difficult but have greater implacability in some special kind of mathematical problems.


\subsection{Indefinite integral}

 For all positive  integer values of n

$ \int {e^{x^{n}}}dx $\\

$= e^{x^n } x - n\int {e^{x^n } .x^n } dx $\\

$= e^{x^n } x - e^{x^n } \frac{n}{{n + 1}}x^{n + 1}  + \frac{{n^2 }}{{n + 1}}\int {e^{x^n } .x^{2n} }dx $ \\

$= e^{x^n } x - e^{x^n } \frac{n} {{n + 1}}x^{n + 1}  + e^{x^n }\frac{{n^2 }} {{(n + 1)(2n + 1)}}x^{2n + 1}  - \frac{{n^3 }} {{(n + 1)(2n + 1)}}\int {e^{x^n } .x^{3n} } dx \hfill $\\

$= e^{x^n } x - e^{x^n } \frac{n} {{n + 1}}x^{n + 1}  + e^{x^n }\frac{{n^2 }} {{(n + 1)(2n + 1)}}x^{2n + 1} $\\

$\hfill - e^{x^n } \frac{{n^3}} {{(n + 1)(2n + 1)(3n + 1)}}x^{3n + 1}+ \frac{{n^4 }}{{(n + 1)(2n + 1)(3n + 1)}}\int {e^{x^n } .x^{4n} } dx $ \\

$= e^{x^n} \left[{\sum\limits_{r = 0}^\infty  \frac{{( - 1)^r n^r x^{1 + nr} }}{{(1 + 0n)(1 + 1n)(1 +2n)...(1 + rn)}} } \right] + C_1 \hfill  $\\

$= e^{x^n } \left[ {\sum\limits_{r = 0}^\infty  {\frac{{( - 1)^r n^r x^{1 + nr} }}{{\prod\limits_{p = 0}^r {(1 + pn)} }}} } \right] + C_1 \hfill {\forall x; \; \;  atleast \; \; \forall n \in \mathbb{Z}^ +  } \;\;\;\;\; (3)$ \\

Similarly:\\

$\int {e^{ - x^n } } dx \hfill $ \\
$= e^{ - x^n } \left[ {x + \sum\limits_{r = 1}^\infty  {\frac{{n^r x^{1 + nr} }}{{(1 + 0n)(1 + 1n)(1 + 2n)...(1 + rn)}}} } \right] + C_2 \hfill $ \\
$= e^{ - x^n } \left[ \sum\limits_{r = 0}^\infty  \frac{n^r x^{1 + nr} } {\prod\limits_{p = 0}^r {(1 + pn)}}  \right] + C_2  \hfill \forall x; \; \;  atleast \; \; \forall n \in \mathbb{Z}^ +    (4)$ \\

Note that equation no. (3) and (4) are valid for all x for all n other than${\left[{ - 1 \leqslant n \leqslant 0} \right]} $

\subsubsection{Observation:}

Taking n= odd and putting x= -y i.e. dx= -dy in equation no. (3) we get,\\

$\int {e^{( - y)^n } } d( - y) \hfill $\\
$= e^{( - y)^n } \left[ {\sum\limits_{r = 0}^\infty  {\frac{{( - 1)^r ( - 1)^{1 + nr} n^r y^{1 + nr} }}
{{\prod\limits_{p = 0}^r {(1 + pn)} }}} } \right] + C_1  \hfill $\\
$= e^{ - y^n } \left[ {\sum\limits_{r = 0}^\infty  {\frac{{( - 1)^{1 + (n + 1)r} n^r y^{1 + nr} }}
{{\prod\limits_{p = 0}^r {(1 + pn)} }}} } \right] + C_1  \hfill $ \\
$=  - e^{ - y^n } \left[ {\sum\limits_{r = 0}^\infty  {\frac{{n^r y^{1 + nr} }} {{\prod\limits_{p = 0}^r {(1 + pn)} }}} } \right] +
C_1   \hfill {\; \; \; as \; n \; is \; assumed \; to \; be \; odd} $\\

or,
$\;\;\;  \int {e^{ - y^n } } dy = e^{ - y^n } \left[ {\sum\limits_{r = 0}^\infty  {\frac{{n^r y^{1 + nr}}}{{\prod\limits_{p = 0}^r {(1 + pn)} }}} } \right] - C_1  \hfill $\\

or,
$\;\;\;  \int {e^{ - x^n } } dx = e^{ - x^n } \left[ {\sum\limits_{r = 0}^\infty  {\frac{{n^r x^{1 + nr} }}
{{\prod\limits_{p = 0}^r {(1 + pn)} }}} } \right] - C_1 \hfill (4A)$\\

Comparing equation no. (4) and (4A) we have $C_{1}$ \; =\; $-C_{2}$

\subsection{Definite integral}

\subsubsection{Proper integrals:}

For proper integrals, convergence of $ e^{x^n } \left[
{\sum\limits_{r = 0}^\infty  {\frac{{( - 1)^r {n^r} x^{1 + nr}
}}{{\prod\limits_{p = 0}^r {(1 + pn)} }}} } \right] $     and    $
e^{x^n } \left[{\sum\limits_{r = 0}^\infty  {\frac{{n^r x^{1 +
nr}}}{{\prod\limits_{p = 0}^r {(1 + pn)} }}} } \right] $  for finite
values of x, are essential and convergence of $ {e^{x^{n} }} {
\left[ {\sum\limits_{r = 0}^{\infty} {\frac{{{( - 1)}^r} {n^{r}}
{x^{1 + nr}}} {\prod\limits_{p = 0}^r {(1 +pn)} }} } \right]} $
depends upon both \quad $ e^{x^n}  $ \quad
and \quad $ \left[ {\sum\limits_{r = 0}^{\infty}  {\frac{{( - 1)^r {n^r} x^{1 + nr} }}{{\prod\limits_{p = 0}^r {(1 + pn)} }}} } \right].$ \\

We thus performed the ratio test of both $ e^{x^n } $ and $ \left[ {\sum\limits_{r = 0}^{\infty}  {\frac{{( - 1)^r n^r x^{1 + nr} }}{{\prod\limits_{p = 0}^r {(1 + pn)} }}} } \right]$

which shows.

$\quad \mathop {\lim }\limits_{r \to \infty } \left| {\frac{{u_{r + 1} }}{{u_r }}} \right| = \mathop {\lim }\limits_{r \to \infty } \left| {\frac{{x^{nr + n} r!}}{{x^{nr} (r + 1)!}}} \right| \quad  = \mathop {\lim }\limits_{r \to \infty } \frac{1}{{r + 1}}\left| {x^n } \right| = 0( < 1) \hfill \forall x $ \\

$\quad = \mathop {\lim }\limits_{r \to \infty } \left| {\frac{{u_{r + 1} }}{{u_r }}} \right| = \mathop {\lim }\limits_{r \to \infty } \left| { - \frac{{n^{r + 1} x^{1 + nr + n} \prod\limits_{p = o}^{r + 1} {(1 + pn)} }}{{n^r x^{1 + nr} \prod\limits_{p = o}^r {(1 + pn)} }}} \right| $

$\quad = \mathop {\lim }\limits_{r \to \infty } \left| { - \frac{{nx^n }}{{1 + n + rn}}} \right| $

$\quad = \mathop {\lim }\limits_{r \to \infty } \frac{n}{{1 + n + rn}}\left| { - x^n } \right| = 0( < 1) \hfill \forall x $\\

So, it converges for all finite values of x.

Hence, $ e^{x^n } \left[ {\sum\limits_{r = 0}^\infty  {\frac{{( -1)^r n^r x^{1 + nr} }}{{\prod\limits_{p = 0}^r {(1 + pn)} }}} }\right] $ i.e. $ \int {e^{x^n } } dx $ converges for all finite
values of x. \\
Similarly it can be shown that  $ \int {e^{-x^n } } dx $  also converges for all finite values of x. Same conclusion can be drawn from integrals of form of (1) or (2) by ratio tests.

\subsubsection{Modification of interval of convergence: improper integral}

 Since the functions $ e^{x^{n}}$  and $ e^{-x^{n}}$ are having no point of singularity, finite integrals are improper when any (or both)  of the limits are infinite.

Now,

$ \int\limits_{ - \infty }^{ + \infty } {e^{x^n } } dx = \int\limits_{ - \infty }^0 {e^{x^n } } dx + \int\limits_0^{ + \infty} {e^{x^n } } dx \hfill (5) $\\

Now, by $\mu$ test  of $ \int\limits_0^\infty  {e^{x^n } } dx $....(6) we see $\mathop {\lim }\limits_{x \to \infty } x^\mu  e^{x^n }  = \infty $  for, $\mu=0$ or 1($\le$1) (where x assumes only positive
value) is evidently diverging.

Again, for $ \int\limits_{ - \infty }^0 {e^{x^n } } dx $ ; x Assumes only negative values in this interval\\
 Let, x=-y; i.e. dx=-dy (y $>$ 0)

\paragraph{for n = even}:

$ \int\limits_{ - \infty }^0 {e^{x^n } } dx = \int\limits_\infty ^0 {e^{y^n } } ( - dy) \;=\; \int\limits_0^\infty  {e^{y^n } } dy = \int\limits_0^\infty  {e^{x^n } } dx...........\left\{ {As,\int\limits_a^b {f(y)dy = \int\limits_a^b {f(x)dx} } } \right\} $

By earlier conclusion it diverges

\paragraph{for n = odd }:

$ \int\limits_{ - \infty }^0 {e^{x^n } } dx = \int\limits_\infty ^0 {e^{ - y^n } } ( - dy) \;=\;  \int\limits_0^\infty  {e^{ - y^n } } dy = \int\limits_0^\infty  {e^{ - x^n } } dx...........\left\{ {As,\int\limits_a^b {f(y)dy = \int\limits_a^b {f(x)dx} } } \right\}$ \\

By $ \mu $ test it is evident that

$ \mathop {\lim }\limits_{x \to \infty } x^\mu  e^{ - x^n } \; = \; \mathop{\lim }\limits_{x \to \infty } \left[ {\frac{{x^2 }}{{e^{x^n } }}} \right] \;\;\; For \;\; n = 2( > 1) \;\;\; = \;\;
\mathop{\lim }\limits_{x \to \infty } \left[ {\frac{1}{{\sum\limits_{r = 0}^\infty  \frac{x^{rn - 2} }{r!} }}} \right] \;\; \hfill (7) $\\

For n $\le$ 0 ,It is $\infty$ ( diverges) (n need not be integer even). But we have not considered so. For n $>$ 0, it can be shown (not necessary that n is odd, it may be even or fraction also) that
expression (7 ) is convergent For integer value of n, the expression obviously is 0, hence convergent. For fractions 'r' varies from 0 to $\infty$ so, what ever the value of n may be, as it is finite, there
will be r's $>$ $\frac{2}{n}$. So,

$ \mathop{\lim }\limits_{x \to \infty } \left[{\frac{1}{{\sum\limits_{r = 0}^\infty  \frac{x^{rn - 2}} {r!} }}} \right] $

$= \mathop {\lim }\limits_{x \to \infty } \left[ \frac{1} {\frac{x^{-2}}{0!} + \frac{x^{n-2}}{1!} + \frac{x^{2n-2}}{2!}+ ...\frac{x^{rn-2}}{r!}+...{\frac{1}{(\frac{2}{n})!}}x^{0} + \infty
+ \infty  + \infty+...} \right] = 0 \hfill  (7A) $\\

Hence convergent. \\
 So, $ \int\limits_{0}^{\infty} {e^{-x^n } } dx$ ,is convergent  $\forall n \in \mathbb{Z}^ +$ i.e. for odd n's also. \\

i.e. $\int\limits_{ - \infty }^0 {e^{x^n } } dx $ converges.

So modified interval of convergence of $ \int\limits_{}^{} {e^{x^n }
} dx $ is $ ]-\infty,\infty[$ for n = even. And
 $- \infty  \leqslant x < \infty$ for n=odd.

Now $ \int\limits_{ - \infty }^{ + \infty } {e^{ - x^n } } dx =
\int\limits_{ - \infty }^0 {e^{ - x^n } } dx + \int\limits_0^{ +
\infty } {e^{ - x^n } } dx \hfill (8) $ \\

As shown earlier, $ \int\limits_0^\infty  {e^{ - x^n } } dx $ is
convergent  $ \forall n \in \mathbb{Z}^ + $ and for,
$\int\limits_{\infty}^{0} {e^{- x^n } } dx $, x Assumes only
negative values in this interval.\\

Let, x=-y;   i.e. dx=-dy  (y$\ge$0)

\paragraph{for n = even}:

$\int\limits_{ - \infty }^0 {e^{ - x^n } } dx = \int\limits_\infty ^0 {e^{ - y^n } } ( - dy)  = \int\limits_0^\infty  {e^{ - y^n } } dy $
$= \int\limits_0^\infty  {e^{ - x^n } } dx \hfill ......... \left\{
{As,\int\limits_a^b {f(y)dy = \int\limits_a^b {f(x)dx} } } \right\}
\hfill $\\

By earlier proof it is convergent.

\paragraph{for n = odd}:

$\int\limits_{ - \infty }^0 {e^{ - x^n } } dx = \int \limits_{\infty} ^0 {e^{y^n } } ( - dy) = \int\limits_0^\infty  {e^{y^n } } dy = \int\limits_0^{\infty}  {e^{x^n } } dx \hfill
......... \left\{ {As,\int\limits_a^b {f(y)dy = \int\limits_a^b
{f(x)dx} } } \right\} \hfill $\\
Which diverges by earlier proof. So modified interval of convergence of $ \int e^{-x^n}  dx $  is $[-\infty,\infty]$ for n= even and $ - \infty  < x \leqslant \infty $ for n=odd.

\subsubsection{Observation}

$ e^{x^n } \left[ {\sum\limits_{r = 0}^\infty  {\frac{{{( - 1)^r} {n^r} x^{(1 + nr)} }}{{\prod\limits_{p = 0}^r {(1 + pn)} }}} } \right] \; =\; \left[\sum\limits_{r = 0}^\infty  {\frac{{x^{nr} }}{{r!}}} \right] \left[ {\sum\limits_{r = 0}^\infty  {\frac{{( - 1)^r n^r x^{1 + nr} }}{{\prod\limits_{p = 0}^r {(1 + pn)} }}} } \right] $\\

 $= \left[ {1 + \frac{{x^n }}{{1!}} + \frac{{(x^n )^2 }}{{2!}} + \frac{{(x^n )^3 }}{{3!}} + ...\infty } \right]\left[ {\frac{{n^0 x^{1 + 0n} }}{{(1 + 0n)}} - \frac{{n^1 x^{1 + 1n} }}{{(1 + 0n)(1 + 1n)}} + \frac{{n^2 x^{1 + 2n} }}{{(1 + 0n)(1 + 1n)(1 + 2n)}}...\infty } \right] $\\
 $= \left[ {\frac{{n^0 }}{{0!(1 + 0n)}}} \right]x^{1 + 0n}  - \left[ {\frac{{n^1 }}{{0!(1 + 0n)(1 + 1n)}} - \frac{{n^0 }}{{1!(1 + 0n)}}} \right]x^{1 + 1n} $ \\
 $\hfill + \left[ {\frac{{n^2 }}{{0!(1 + 0n)(1 + 1n)(1 + 2n)}} - \frac{{n^1 }}{{1!(1 + 0n)(1 + 1n)}} + \frac{{n^0 }}{{2!(1 + 0n)}}} \right]x^{1 + 2n}  \;\;\;   .....\infty $\\
 $= \sum\limits_{m = 0}^\infty  {\left( { - 1} \right)^m \left[ {\sum\limits_{r = 0}^m {\frac{{\left( { - 1} \right)^r n^{m - r} }}{{r!\prod\limits_{p = 0}^{m-r} {(1 + pn)} }}} } \right]} x^{1 + mn} \hfill (9) $\\

Similarly \\

$ e^{x^n } \left[{\sum\limits_{r=0}^\infty {\frac{n^r x^{1 +
nr}}{\prod\limits_{p = 0}^r {(1 + pn)}}}}\right] \; =\;
\sum\limits_{m = 0}^\infty  {\left[{\sum\limits_{r=0}^m
{\frac{\left({-1}\right)^r n^{m-r}}{r! \prod\limits_{p=0}^{m-r}
{(1+pn)}}}} \right]} x^{1+mn} \hfill (10)$\\

From equation no.(1)and(9), comparing coefficients of $ x^{1 + mn}$
we get,\\

$\frac{1}{m!(nm+1)}=\left({-1}\right)^m \left[{\sum\limits_{r=0}^m
{\frac{\left({-1} \right)^r n^{m-r}}{r!\prod\limits_{p=0}^{m-r}
{(1+pn)}}}} \right] \hfill  (11) $\\

 From equation no. (2)and(10) comparing coefficients of $ x^{1 +
mn}$ we get, \\

$ \left({-1}\right)^m \frac{1} {{m!(nm + 1)}}  = \left[
{\sum\limits_{r = 0}^m {\frac{\left( { - 1} \right)^r n^{m - r} }
{r!\prod\limits_{p = 0}^{m-r} {(1 + pn)} }} } \right] \hfill (12) $\\

Obviously equation (11) gives equation (12) if multiplied by $ \left( { - 1} \right)^m $. \\

\section{Application of integrals of $e^{x^{n}}$ and $e^{-x^{n}}$ in the form of series as series solution of some differential equations.}

Considering two second order differential equations:\\

$x\frac{{d^2 y}}{{dx^2 }} - (n - 1) \frac{{dy}}{{dx}} - n^2 x^{2n -
1} y - nx^n  = 0 \hfill $

or,

$\frac{{d^2 y}}{{dx^2 }} - \frac{{(n - 1)}}{x}\frac{{dy}}{{dx}} -n^2
x^{2n - 2} y = nx^{n - 1}  \hfill (13) $

and

$x\frac{{d^2 y}}{{dx^2 }} - (n - 1)\frac{{dy}}{{dx}} - n^2 x^{2n - 1} y + (n - 1)= 0 \hfill $

or,

$\frac{d^2 y}{dx^2 } - \frac{(n - 1)}{x}\frac{dy}{dx} - n^2 x^{2n - 2} y = \frac{ - (n - 1)}{x} \hfill (14)$

Here equation no. (13), (14) are actually the special forms of the equation\\   $\frac{d^2 y}{dx^2 }+ P\frac{dy}{dx}+ Qy =R \;\;\;\;\;\;\;\;\;\; (15)$ Where P,Q,R are functions of x only.

By method of change of independent variable and taking some arbitrary function $ z= \int {e^{ - \int {Pdx} } } dx$    equation no. (15) can be reduced to $ \frac{d^2 y}{dz^2 } + Q_1 y = R_1 \hfill (15A) $

where,    $ {\rm{Q}}_{\rm{1}}  = \frac{Q}{\left( {\frac{dz}{dx}} \right)^2 }, \;\;\;\; R_1  = \frac{Q}{\left(
{\frac{dz}{dx}} \right)^2 } $

Now, for both (13) and (14) $ z = \frac {x^{n}}{n} $ i.e. $x = \left({nz} \right)^{\frac{1}{n}} \; =\;$ let,$\;\; f^{ - 1} (z). $

Now, $\left( {\frac{{dz}}{{dx}}} \right)^2 = x^{2n - 2}$  or, $Q_1 =  - n^2 $

So, (13) becomes  $\frac{{d^2 y}}{{dz^2 }} - y = R_1 \hfill (16) $

Now $ R_1  = nx^{1 - n}  = \frac{n}{{\left[ {f^{ - 1} (z)} \right]^{n - 1} }} \;\;\;$ for equation (13) and =$ \frac{{ - (n - 1)}}{{\left[ {f^{ - 1} (z)} \right]^{2n - 1} }} \;\;\;$ for equation (14)

So, equation \;\; (13)\; becomes

$ \frac{{d^2 y}}{{dz^2 }} - n^2 y = \frac{n}{{\left[ {f^{ - 1} (z)}
\right]^{n - 1} }} \hfill (17) $

and equation   (14) becomes

$ \frac{{d^2 y}}{{dz^2 }} - n^2 y = \frac{{ - (n - 1)}}{{\left[ {f^{
- 1} (z)} \right]^{2n - 1} }} \hfill (18) $

Evidently (17) and (18) are second order linear differential
equation with constant coefficients. The complete solutions of these
are given by Complimentary function(C.F.)+Particular integral(P.I.)

For both, C.F. is solution of the equations.

$ \frac{{d^2 y}}{{dz^2 }} - n^2 y = 0 $  or, symbolically,  $(D^2  - n^2 )y = 0 \hfill \left\{ {D = \frac{d}{{dz}}} \right\} $

or,
C.F.$ = k_1 e^{nz}  + k_2 e^{ - nz}  = k_1 e^{x^n }  + k_2 e^{ - x^n }$        \hfill       $k_1$, $k_2$, arbitrary, constants. \\

Now P.I. of equation no. ...(17) is  $(P.I.)_1$ = $\frac{1}{{(D^2  - n^2 )}}\left[ {\frac{n}{{\left\{ {f^{ - 1} (z)} \right\}^{n - 1} }}} \right]$ \\
 and for equation no. ...(18),  $(P.I.)_2$ = $ \frac{1}{{(D^2  - n^2 )}}\left[ {\frac{{ - (n - 1)}}{{\left\{ {f^{ - 1} (z)} \right\}^{2n - 1} }}} \right] $ \\

It seems to be complicated to be evaluated in this approach. So alternative methods can be taken. By definition particular integral is any particular solution of the given differential equation. So
$(P.I.)_1$ is a particular solution of  (17) or (13), where as $(P.I.)_2$  is any particular solution of equation no. (18) or (14).

It can be shown that $ g(x) = \left[ {\sum\limits_{r = 0,2,4,..\infty } {\frac{{n^r x^{1 + nr} }}{{\prod\limits_{p = 0}^r {(1 + pn)} }}} } \right] $ is a particular solution of equation no.
(13).

Where as, $ f(x) = \left[ {\sum\limits_{r = 1,3,5,..\infty } {\frac{{n^r x^{1+ nr} }}{{\prod\limits_{p = 0}^r {(1 + pn)} }}} } \right]$ is a particular solution of equation (14).

By equation no. (3) and (4)

$\int {e^{x^n } } dx = e^{x^n } \left[ {f(x) - g(x)} \right] + C_1
\hfill (19) $

and

$\int {e^{ - x^n } } dx = e^{ - x^n } \left[ {f(x) + g(x)} \right] + C_2 \hfill (20) $

Differentiating equation (19) with respect to x,

$e^{x^n }  = e^{x^n } \left[ {f'(x) - g'(x)} \right] + e^{x^n } \left[ {f(x) - g(x)} \right]nx^{n - 1} \hfill (21) $

Differentiating equation (20) with respect to x,

$e^{ - x^n }  = e^{ - x^n } \left[ {f'(x) + g'(x)} \right] - e^{ - x^n } \left[ {f(x) + g(x)} \right]nx^{n - 1} \hfill (22) $

or

$\left[ {f'(x) - g'(x)} \right] + \left[ {f(x) - g(x)} \right]nx^{n - 1}  = 1 \hfill (23) $

and

$\left[ {f'(x) + g'(x)} \right] - \left[ {f(x) + g(x)} \right]nx^{n - 1}  = 1 \hfill (24) $

from equation (23) $\&$ (24)

$f'(x) = nx^{n - 1} g(x) + 1 \hfill (25) $

$g'(x) = nx^{n - 1} f(x) \hfill (26) $

Differentiating equation (25) $\&$ (26) and using (25) and (26) we
have,

$ g''(x) - \frac{{(n - 1)}}{x}g'(x) - n^2 x^{2n - 2} g(x) = nx^{n - 1} \hfill (27) $

and

$ f''(x) - \frac{{(n - 1)}}{x}f'(x) - n^2 x^{2n - 2} f(x) = \frac{{ - (n - 1)}}{x} \hfill (28) $\\

So,

$ g(x) = \left[ {\sum\limits_{r = 0,2,4,..\infty } {\frac{{n^r x^{1 + nr} }}{{\prod\limits_{p = 0}^r {(1 + pn)} }}}} \right]$ is a particular solution of (13) \\

where as,

$ f(x) = \left[ {\sum\limits_{r = 1,3,5,..\infty } {\frac{{n^r x^{1 + nr} }}{{\prod\limits_{p = 0}^r {(1 + pn)} }}} } \right]$ is a particular solution of (14).\\

Now the complete solution of (13) is  $ y = k_1 e^{x^n}  + k_2 e^{ -x^n}  + \left[ {\sum\limits_{r = 0,2,4,..\infty } {\frac{{n^r x^{1 + nr} }}{{\prod\limits_{p = 0}^r {(1 + pn)} }}} } \right] \hfill (28A)$ \\

And that of (14) is $ y = k_1 e^{x^n}  + k_2 e^{ -x^n}  + \left[ {\sum\limits_{r = 1,3,5,..\infty } {\frac{{n^r x^{1 + nr} }}{{\prod\limits_{p = 0}^r {(1 + pn)} }}} } \right]\hfill (28B)$\\

Where, $k_1$ and $k_2$ are arbitrary constants.

\subsection{Graphical representation of   y(x) = $e^{ - x^n }$ for even values of n.}

\begin{figure}[!h]
\centering
\includegraphics[width=9 cm]{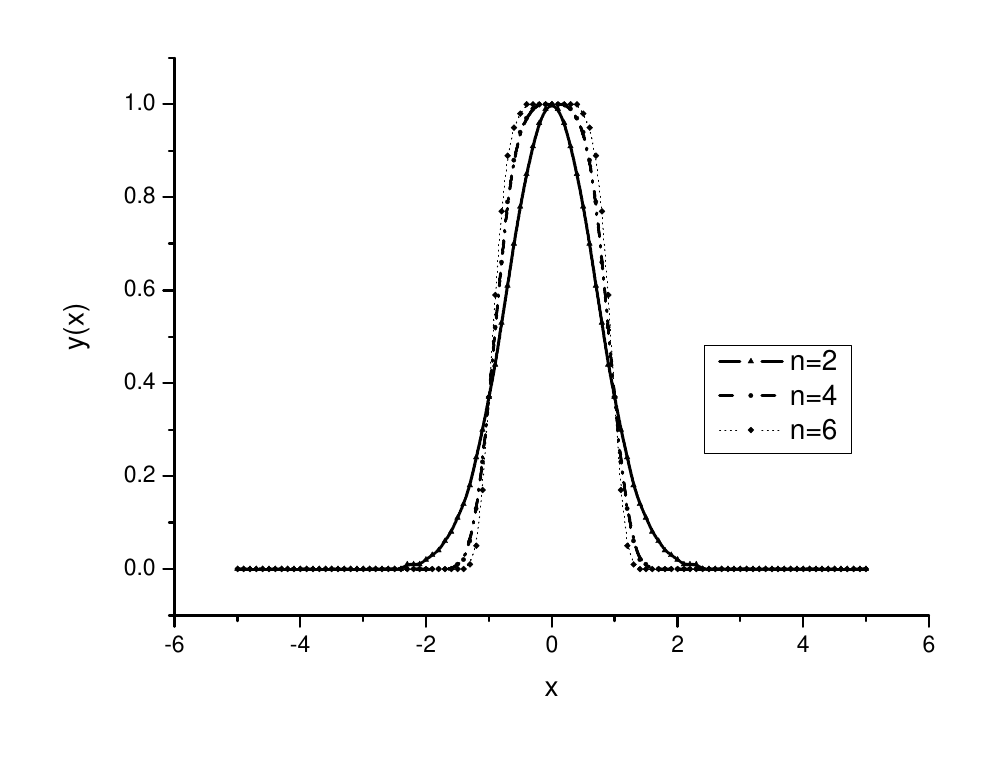}
\caption{Values of y(x) for n= 2,4,6}
\end{figure}

Fig 1. shows the curves of y(x) =  for different even values of n.

Again, we have,

When n tends to $\infty$, in that case

\begin{itemize}

\item  $ y(x) = \mathop {\lim }\limits_{n \to \infty } e^{ - x^n }  = 1$, for - 1$ < x <$ 1

\item   = $\frac{1}{e}$,  at,x = $\pm 1 $

\item   = 0,  elsewhere

\end{itemize}

Which describes a rectangular function [fig 2.].

\begin{figure}[!h]
\centering
\includegraphics[width=8 cm]{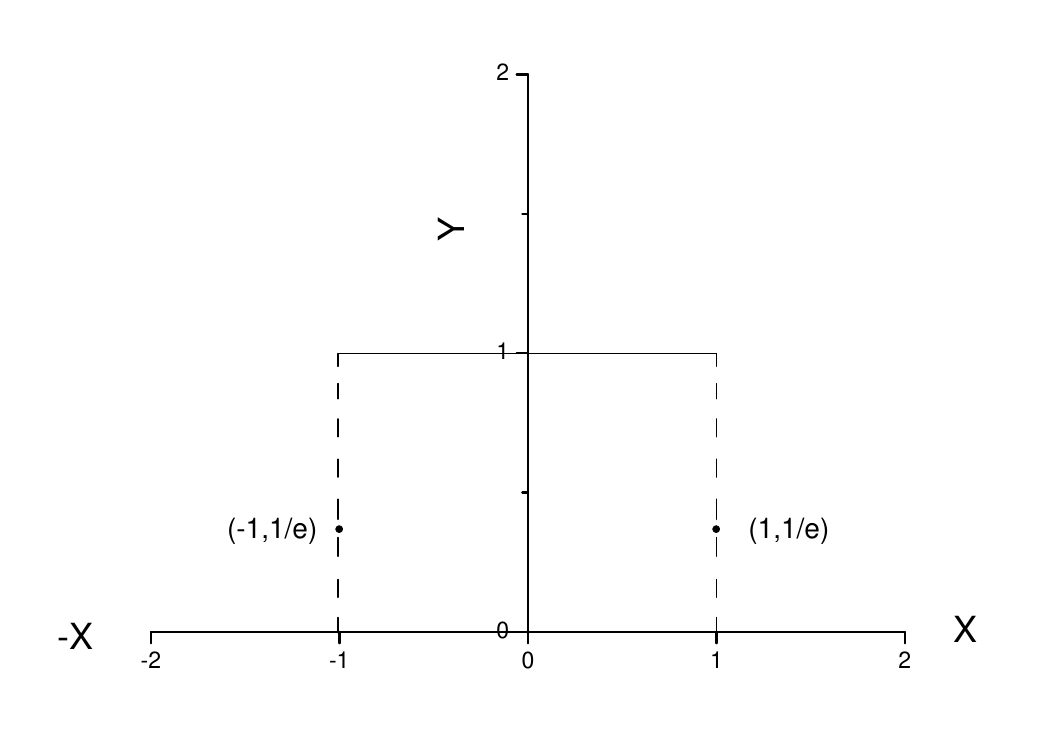}
\caption{ $ y(x) = \mathop {\lim }\limits_{n \to \infty } e^{ - x^n } $}
\end{figure}

But if n is very big(yet finite),

$ \mathop {\lim }\limits_{x \to  \pm 1 + } e^{ - x^n }  = \mathop {\lim }\limits_{x \to  \pm 1 - } e^{ - x^n }  = \mathop {\lim}\limits_{x \to  \pm 1  } e^{ - x^n }  = \frac{1}{e}$ \\

To be noted that taking a big value of even n the aforesaid rectangular function can be approximately replicated in the field of computation. The value of n is to be chosen depending upon the ability of the computing machine.

\section{Admissibility of $ y(x) = Ae^{ - \frac{1}{n} {\left( {\frac{x - m}{\sigma }} \right)}^n}$ as a Probability Density Function (P.D.F) in $(-\infty <x< \infty) $ for positive, even values of n.}

Let us, consider a function $ y(x) = A {e^{- \frac{1}{n}\left( \frac{x - m}{\sigma } \right)}}^{n}  $ (A being constant) y(x) must hold the following properties in order to be a Probability Density Function (P.D.F) in (a,b).\\
 I. y(x)$>$0 in(a,b),  II. y(x)=0 else where and III. $ \int\limits_{ - \infty }^\infty  {y(x)} dx = 1 $ \\

It is known that normal distribution is a bell shaped continuous probability distribution and it must hold the following properties.

\begin{itemize}

\item Normal curve has two parameters, mean = median = mode = 'm' and standard deviation $\sigma$.

\item normal curve is of form $ y(x) = \frac{1}{\sigma \sqrt {2\pi }}e^{ - \frac{1}{2}\left( {\frac{x - m}{\sigma }} \right)^2 } \hfill (29).$

\item the curve is symmetrical about x = m.

\item odd order central moments are 0.

\item Coefficient of Skewness and Kurtosis excess are 0.
\end{itemize}

Here, the Probability density function in form of standardized variable $z = \frac{{x - m}}{\sigma }$ is $ y(z) = \frac{1}{{\sqrt {2\pi } }}e^{ - \frac{z}{2}^2 } \hfill (30).$

For $ y(x) = Ae^{ - \frac{1}{n}\left( {\frac{{x - m}}{\sigma }}\right)^n } $ , (n being even) conditions I and II are automatically satisfied. In order to satisfy condition III, constant A is to be
chosen suitably. So that,

$\int_{ - \infty }^\infty  {y(x)} dx = 1 $

or,

$ \int\limits_{ - \infty }^\infty  {Ae^{ - \frac{1}{n}\left( {\frac{{x - m}}{\sigma }} \right)^n } } dx = 1 $

or,

$ A\sqrt[n]{n}\sigma \int\limits_{ - \infty }^\infty  {e^{ - \left({\frac{x - m}{\sqrt[n]{n}\sigma }} \right)^n } } d\left( {\frac{x - m}{\sqrt[n]{n}\sigma }} \right) = 1 $

For even values of n, it has been shown that the improper integral $ \int\limits_{ - \infty }^{ + \infty } {e^{ - x^n } } dx$   exists and that can be expressed as 2 $\Gamma$ (1 + $ \frac{1}{n}$).... (30A). Let that be denoted as $P_n$.

So,

$ A\sqrt[n]{n}\sigma \int\limits_{ - \infty }^\infty  {e^{ - \left( {\frac{{x - m}}{{\sqrt[n]{n}\sigma }}} \right)^n } } d\left( {\frac{{x - m}}{{\sqrt[n]{n}\sigma }}} \right) = A\sqrt[n]{n}\sigma P_n  $

 or, $ A\sqrt[n]{n}\sigma P_n  = 1 $  or,  $ A = \frac{1}{{\sqrt[n]{n}\sigma P_n }} $

Therefore,

$ y(x) = \frac{1}{{\sqrt[n]{n}\sigma P_n }}e^{ - \frac{1}{n}\left( {\frac{{x - m}}{\sigma }} \right)^n }$  is a Probability density function in $(-\infty <x< \infty ) \hfill (31) $

 Which can be called $j^{th}$, order (where n=2j) Normal Distribution. Equation (29) is a $1^{st}$ order Normal Distribution in this regard.

\subsection{ Properties of distribution corresponding to P.D.F $ y(x)= \frac{1}{\sqrt[n]{n}\sigma P_n }e^{ - \frac{1}{n}\left( {\frac{x - m}{\sigma }} \right)^n }$}

Equation no. ..(31) can be further reduced in  form of standardized variable z = $\frac{{x - m}}{\sigma }$ is    $y(z) = \left[ \frac {1}{{\sqrt[n]{n}P_n }}\right]  e^{ - \frac{z}{n}^n } \;\;\; (-\infty
<z< \infty) \hfill (32).$

\subsubsection{The distribution is symmetric about x=m i.e. z=0 and is bell shaped.}

From equation (32) we have,

$ y(z) = \frac{1}{{\sqrt[n]{n}P_n }}e^{ - \frac{z}{n}^n } \hfill ( - \infty  < z < \infty ) $

For even values of n, this is an even function. So, this is symmetric about z=0. Which implies equation (32) is symmetric about x=m.

Again, we have,

$ y(z) = \frac{1}{\sqrt[n]{n}P_n }e^{ - \frac{z}{n}^n }  $

or,  $ \frac{{dy(z)}}{{dz}} \; = \; \frac{1}{{\sqrt[n]{n}P_n }}e^{ -\frac{z}{n}^n } \left( { - z^{n - 1} } \right) \; =\;  - z^{n - 1} y $

So, $ \frac{{dy(z)}}{{dz}}= 0$,   at [z = 0,and,y = 0]

at  z = 0,    $ y(z) = \frac{1}{{\sqrt[n]{n}P_n }} $  and  at y = 0, $ e^{ - \frac{z}{n}^n }  = 0 $  or, $ z =\pm \infty \hfill (n= even) $

i.e.  $ \frac{dy(z)}{dz} = 0, $   at   $ \left[\left( {0,\frac{1}{\sqrt[n]{n}P_n }} \right)\right]$,  and,   $ \left[\left( { \pm \infty ,0} \right) \right] $

Now,

$\frac{{d^2 y(z)}}{{dz^2 }} = \frac{{d\left( { - z^{n - 1} y} \right)}}{{dz}} $

$=  - \left[ {z^{n - 1} \frac{{dy}}{{dz}} + y(n - 1)z^{n - 2} } \right] $

$=  - \left[ {z^{n - 1} \left( { - z^{n - 1} y} \right) + y(n - 1)z^{n - 2} } \right] $

$=  z^{n - 2} \left( {z^n -(n - 1)} \right)y \hfill (33) $

Obviously, this expression approaches to 0 from negative side of number line.

Let $\epsilon$ be a variable, which tends to z or,  $ \left[ {\frac{d^2 y(z)}{dz^2 }} \right]_{\left(
{0,\frac{1}{\sqrt[n]{n}P_n }} \right)}  = \mathop {\lim }\limits_{\epsilon  \to 0} \frac{1}{\sqrt[n]{n}P_n }\left[ {\epsilon ^{2n - 2}  - \left( {n - 1} \right)\epsilon ^{n - 2} }\right] $

Now, n  $ \ge 2,or,(n - 1) \ge 1 $ , 2n-2 $>$ n-2 and $\epsilon$ is a fraction which tends to zero.

$ \mathop {\lim }\limits_{\epsilon  \to 0} \frac{1}{\sqrt[n]{n}P_n }\left[ {\epsilon ^{2n - 2}  - \left( {n - 1} \right)\epsilon ^{n - 2} } \right] < 0 $

or,$\left[ {\frac{d^2 y(z)}{dz^2 }} \right]_{\left( {0,\frac{1}{\sqrt[n]{n}P_n }} \right)}  < 0 $

or,  y is maximum at z=0 and is  $\frac{1}{\sqrt[n]{n}P_n } . $\\

Again,

Let '$\delta$' be a variable, which tends to z and $\gamma$ is a variable tending to y.

or, $\left[ {\frac{d^2 y(z)}{dz^2 }} \right]_{\left( { \pm \infty ,0} \right)}  = \mathop {\lim }\limits_{\delta  \to  \pm \infty } \mathop {\lim }\limits_{\gamma  \to 0} \gamma \delta ^{n - 2} \left[ {\delta ^n   - \left( {n - 1} \right) } \right] $\\

Now, n is even number and $\ge$ 2  and '$\delta$' is a number which tends to $\pm \infty$.\\

or, $ \left[ {\frac{d^2 y(z)}{dz^2 }} \right]_{\left( { \pm \infty ,0} \right)} \; =\; \mathop {\lim }\limits_{\delta  \to  \pm \infty } \mathop {\lim }\limits_{\gamma  \to 0} \gamma \left[{\delta ^{2n - 2}  - \left( {n - 1} \right)\delta ^{n - 2} } \right]  \hfill (34) $\\

Obviously expression (34) approaches 0 from positive side of number line, as $\gamma$ is non negative as it represents y (evident from the expression (32)).

or, $ \left[ {\frac{{d^2 y(z)}}{{dz^2 }}} \right]_{\left( { \pm \infty ,0} \right)}  > 0 $

So, y is minima at z= $\pm \infty$ and is 0.

Therefore, at  $ \left( {0,\frac{1}{{\sqrt[n]{n}P_n }}} \right)$ the function gets maxima and at $ \left( { \pm \infty ,0} \right)$ the function gets minima. \\
Now equation (33)=0 at points of inflexions. The condition is satisfied for z=0, y=0. or, $ \left[ {z^n  - \left( {n - 1} \right)} \right] = 0 $. At first two conditions, maxima and minima lie. So condition for point inflexion is  $z^n  - \left( {n - 1} \right)=0 $  or, $\left( {z^2 } \right)^{\frac{n}{2}}  - \left( {n - 1} \right) = 0 \hfill (34A).$\\

$\frac{n}{2}$ is whole number as n is even. So, equation (34A) may have $\frac{n}{2}$ numbers of possible solutions, but at least one possible real solution is.

$\left( {z^2 } \right) = \left( {n - 1} \right)^{\frac{2}{n}}$  or,$z =  \pm \left( {n - 1} \right)^{\frac{1}{n}} $\\

there, $ y = e^{ - \frac{1}{n}\left\{ { \pm \left( {n - 1} \right)^{\frac{1}{n}} } \right\}^n } \; = \; e^{ - \frac{(n - 1)}{n}} \; = \; \frac{1} {e}e^{\frac{1}{n}}  \hfill $\\

So, the curve has two points of inflexions, one in first quadrant $\left( {\left( {n - 1} \right)^{\frac{1}{n}}} ,\frac{1} {e}e^{\frac{1}{n} } \right)$ another is in second quadrant $\left( {
- \left( {n - 1} \right)^{\frac{1}{n}}} ,\frac{1} {e}e^{\frac{1}{n}} \right)$. (Where,n is even no).

\begin{figure*}[!h]
\centering
\includegraphics[width=9 cm]{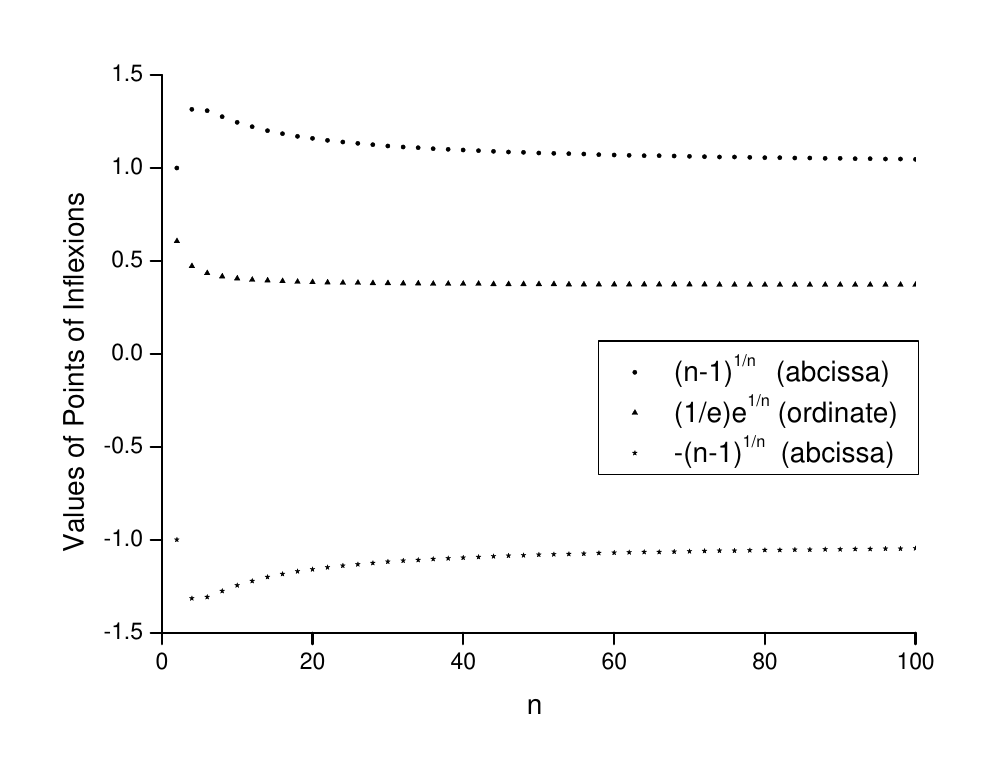}
\caption{Abscissas and ordinates for different values of n of `points of inflexion'.}
\end{figure*}

\footnote{All the graphs are plotted using Origin 7 software package.}

Fig.3 shows variations of abscissas and ordinates of 'points of inflexion' for different values of n.

Again,

$ \frac{{dy(z)}}{{dz}} =  - z^{n - 1} y $, n-1 is an odd no; z varies from $-\infty$ to $\infty$ and y is non-negative,[evident from expression (32)]. So the gradient of the curve of P.D.F [expression (32)] is positive when z is negative and is negative when z is positive.

With the aforesaid arguments it is clear that the distribution is symmetric about x=m i.e. z=0 and is bell shaped.

\subsubsection{Moments about mean.}

\begin{itemize}

\item Odd order moments about mean.

\end{itemize}

$ m_{2p + 1} =  \frac{1}{{p_n \sqrt[n]{n}\sigma }}\int\limits_{ -\infty }^\infty  {(x - m)^{2p + 1} } e^{ - \frac{1}{n}\left( {\frac{{x - m}}{\sigma }} \right)^n } dx \hfill $

$= \frac{1}{{p_n \sqrt[n]{n}\sigma }} \int\limits_{ - \infty }^\infty  {\sigma ^{2p + 1} z^{2p + 1} } e^{ - \frac{z}{n}^n } \sigma dz \hfill $

Integrand is an odd function. So, all odd order moments about mean are 0.

\begin{itemize}

\item Even order moments about mean,

\end{itemize}

We worked out analytically the expressions of even order moments as follows..

$ m_{2p}  = \frac{1} {{p_n \sqrt[n]{n}\sigma }}\int\limits_{ - \infty }^\infty  {(x -m)^{2p} } e^{ - \frac{1} {n}\left( {\frac{{x - m}} {\sigma }} \right)^n } dx  $

$= \frac{1} {{p_n \sqrt[n]{n}\sigma }}\int\limits_{ - \infty }^\infty  {\sigma ^{2p} z^{2p} } e^{ - \frac{{z^n }} {n}} \sigma dz \hfill ...........\left\{ {\frac{{x - m}} {\sigma } = z}\right\}$

$= \frac{{\sigma ^{2p} }} {{p_n \sqrt[n]{n}}}\int\limits_{ - \infty }^\infty  {z^{2p} } e^{ - \frac{{z^n }} {n}} dz \hfill(35A)$

As both `n' and `2p' are even so integrand becomes even function. So (35A) can alternatively written as..

$ = \frac{{2\sigma ^{2p} }} {{p_n \sqrt[n]{n}}}\int\limits_0^\infty {z^{2p}} e^{- \frac{{z^n }} {n}} dz \hfill (35B) $

The integral $ \int\limits_0^\infty  {z^{2p} } e^{ - \frac{{z^n }} {n}} dz $ exists for any value of `p'.

Because, (by $\mu$ test) $  \mathop {\lim }\limits_{z \to \infty }z^\mu z^{2p} e^{ - \frac{{z^n }} {n}} $ tends to zero for �$\mu$=2�($>$1).So integral as stated in (35A) or (35B) exists for any value of `p'.

Now (35A) can be written as..

$ \frac{{\sigma ^{2p} }} {{p_n \sqrt[n]{n}}}\left[ {\left| { - z^{2p - \overline {n - 1} } e^{ - \frac{{z^n }} {n}} } \right|_{ - \infty }^\infty   + \left({2p - \overline {n - 1} } \right)\int\limits_{ - \infty }^\infty {z^{2p - n} } e^{ - \frac{{z^n }}{n}} dz} \right] \hfill  $

$   = \frac{{\sigma ^n }}{{p_n \sqrt[n]{n}}}\left( {2p - \overline {n - 1} } \right)\int\limits_{ - \infty }^\infty  {\sigma ^{2p - n} z^{2p - n}} e^{ - \frac{{z^n }} {n}} dz \hfill $

$   = \left( {\sigma ^n } \right)\left( {2p + 1 - n} \right)m_{2p - n}  \hfill  $\\

By repeated reduction..

$ m_{2p}  = \left( {\sigma ^n } \right)^2 \left( {2p + 1 - n} \right)\left( {2p + 1 - 2n} \right)m_{2p - 2n}  \hfill $

$   = \left( {\sigma ^n } \right)^3 \left( {2p + 1 - n} \right)\left( {2p + 1 - 2n} \right)\left( {2p + 1 - 3n} \right)m_{2p - 3n}  \hfill $

$   = \left( {\sigma ^n } \right)^4 \left( {2p + 1 - n} \right)\left( {2p + 1 - 2n} \right)\left( {2p + 1 - 3n} \right)\left( {2p + 1 - 4n} \right)m_{2p - 4n} ..etc  \hfill(36)$

For some positive integer value of `k' if `2p-kn' assumes zero value, then only by repeated reduction the expression (36) gets reduced to

$ m_{2p}  = \left( {\sigma ^n } \right)^{2p/n} \prod\limits_{k =1}^{2p/n} {\left( {2p + 1 - kn} \right)m_0 } \hfill (36A)$

Obviously, $ \left[ m_0  = 1 \right] $

If there be no such integer `k' for which `2p-kn' becomes `zero', then following things happen. Let `q' be the greatest integer for which `2p-qn $>$ 0'. `2p-qn' must be lesser than `n', else another `n' could be subtracted. As `p, q, n' all are integers and `n' is even so `2p-qn' is an even integer.

In this case the ultimately reduced expression

$m_{2p}  = \left( {\sigma ^n } \right)^q \prod\limits_{k = 1}^q {\left( {2p + 1 - kn} \right)m_{2p - qn} } \hfill (36B) $ \\

For a given `n' and `p' , 2p-qn can assume any even value at most up to n-2. So expression (36B) involves terms �$m_0$(for qn=2p), $ m_2, m_4, m_6 .......m_{n - 2} $

Obviously these can't be further reduced. For a given `n' , any even moment can be expressed in term of one of these as an ultimately reduced form. So these may be called �fundamental moments". For a given even number `n', there should be `n/2' numbers of such fundamental even order moments. In that respect `$m_0$'is the �first fundamental moment" with value equal to unity. If for a given `n' and `2p' , `$m_0$'is achieved by successive reduction then by (36A),

$ m_{2p}  = \left( {\sigma ^n } \right)^{2p/n} \prod\limits_{k =1}^{2p/n} {\left( {2p + 1 - kn} \right)} \hfill (36C) $

To be noted that, as the order of target even order moment is increased, any particular fundamental moment appears periodically with a period of `n'. For a given `n' ,

$m_0$ is achievable in $\left[ m_n ,m_{2n} ,m_{3n} ,m_{4n} .... \right]$ etc.

$m_2$ is achievable in $\left[ m_{2 + n} ,m_{2 + 2n} ,m_{2 + 3n}, m_{2 + 4n} .... \right]$ etc.

$m_4$ is achievable in $\left[ m_{4 + n} ,m_{4 + 2n} ,m_{4 + 3n}, m_{4 + 4n} ....\right]$ etc.

similarly,

$\left[ m_{n - 2}\right]$ is achievable in $\left[m_{\overline {n - 2} + n}, m_{\overline {n - 2}  + 2n}, m_{\overline {n - 2}  + 3n}, m_{\overline {n - 2}  + 4n} ....\right]$

If for different values of `n', $ m_n ,m_{2n} ,m_{3n} ,m_{4n} ....$ etc (i.e. moments involving ``first fundamental moment'') are evaluated by reduction relation, then they appear in a regular pattern as shown below..

$  n = 2, \hfill $
$  m_2  = (\sigma ^2 )^1 .1;,m_4  = (\sigma ^2 )^2 .1.3;,m_6  = (\sigma ^2 )^3 .1.3.5.etc. \hfill $

$  n = 4, \hfill $
$  m_4  = (\sigma ^4 )^1 .1;,m_8  = (\sigma ^4 )^2 .1.5;,m_{12}  = (\sigma ^4 )^3 .1.5.9.etc. \hfill $

$  n = 6, \hfill $
$ m_6  = (\sigma ^6 )^1 .1;,m_{12}  = (\sigma ^6 )^2 .1.7;,m_{18}  = (\sigma ^2 )^3 .1.7.13.etc \hfill $

Thus for a given `n',

$  m_n  = (\sigma ^n )^1 .1;,m_{2n}  = (\sigma ^n )^2 .1.(1 + n);,m_{3n}  = (\sigma ^n )^3 .1.(1 + n).(1 + 2n)$
Thus, \\
$ m_{kn}  = (\sigma ^n )^k \prod\limits_{r = 0}^{k - 1} {(1 + rn)} \hfill (36D) $

Even order moments can also be expressed in following forms..

By (35B)

\[\begin{gathered}
  m_{2p}  = \frac{{2\sigma ^{2p} }}
{{p_n \sqrt[n]{n}}}\int\limits_0^\infty  {z^{2p} } e^{ -
\frac{{z^n}}{n}} dz \hfill \\
   = \frac{{2\sigma ^{2p} \left( {\sqrt[n]{n}} \right)^{2p} }}
{{p_n }}\int\limits_0^\infty  {\left( {\frac{z} {{\sqrt[n]{n}}}}
\right)^{2p} } e^{ - \left( {\frac{z} {{\sqrt[n]{n}}}} \right)^n }
d\left( {\frac{z}{{\sqrt[n]{n}}}} \right) \hfill \\
   = \frac{{2\sigma ^{2p} \left( {\sqrt[n]{n}} \right)^{2p} }}
{{p_n }}\int\limits_0^\infty  {e^{ - \left( {\frac{z}
{{\sqrt[n]{n}}}} \right)^n } \frac{{\left( {\frac{z}
{{\sqrt[n]{n}}}} \right)^{2p} }} {{n\left( {\frac{z}
{{\sqrt[n]{n}}}} \right)^{n - 1} }}} d\left( {\frac{z}
{{\sqrt[n]{n}}}} \right)^n  \hfill \\
   = \frac{{2\sigma ^{2p} \left( {\sqrt[n]{n}} \right)^{2p - n} }}
{{p_n }}\int\limits_0^\infty  {e^{ - \left( {\frac{z}
{{\sqrt[n]{n}}}} \right)^n } } \left( {\frac{z} {{\sqrt[n]{n}}}}
\right)^{2p + 1 - n} d\left( {\frac{z}{{\sqrt[n]{n}}}} \right)^n  \hfill \\
   = \frac{{2\sigma ^{2p} \left( {\sqrt[n]{n}} \right)^{2p - n} }}
{{p_n }}\int\limits_0^\infty  {e^{ - \left( {\frac{z}
{{\sqrt[n]{n}}}} \right)^n } } \left\{ {\left( {\frac{z}
{{\sqrt[n]{n}}}} \right)^n } \right\}^{\frac{1} {n} + \frac{{2p}}
{n} - 1} d\left( {\frac{z}{{\sqrt[n]{n}}}} \right)^n  \hfill \\
   = \frac{{2\sigma ^{2p} \left( {\sqrt[n]{n}} \right)^{2p - n} }}
{{p_n }}\Gamma \left( {\frac{1} {n} + \frac{{2p}} {n}}
\right).............\left[ {As,\frac{1} {n} + \frac{{2p}}
{n} > 0} \right] \;\;\;\;\;\;\;\;\;\;\;\; \hfill (37)  \\
\end{gathered}
\]

Now from (30A) $ p_n  = \frac{2} {n}\Gamma \left( {\frac{1} {n}}
\right) $

$ m_{2p}  = \left( {\sqrt[n]{n}\sigma } \right)^{2p} \frac{{\Gamma
\left( {\frac{1} {n} + \frac{{2p}} {n}} \right)}} {{\Gamma \left(
{\frac{1} {n}} \right)}} \hfill (38)$

\subsection {Skewness and Kurtosis.}

Skewness is the measure of degree of asymmetry and kurtosis is the measure of degree of peakedness of any distribution. i.e. these two are two parameters telling about the geometry of a distribution. For Normal distribution we know that Coefficient of

$ {\text{Skewness = }}\sqrt {\frac{{m_3^2 }} {{m_2^3 }}}$ ....[1]

$ {\text{kurtosis = }}\frac{{m_4 }} {{m_2^2 }} $ and ${\text{kurtosis excess  = }}\frac{{m_4 }} {{m_2^2 }} - 3$. ....[3]

However Pearson expressed Kurtosis as $ \frac{{m_4 }} {{m_2^2 }} - 3 $. ...[2]

Now, Normal distribution is special form of $ y(x) = \frac{1} {{p_n \sqrt[n]{n}\sigma }}(x - m)^{2p} e^{ - \frac{1} {n}\left( {\frac{{x - m}} {\sigma }} \right)^n }$ when `n=2'. For that Skewness is zero as third order (odd order) central moment is zero. In that case kurtosis is `3' and kurtosis excess is 0. By (38) kurtosis could be evaluates as,

\[\begin{gathered}
  \frac{{\left( {\sqrt[n]{n}\sigma } \right)^4 \frac{{\Gamma \left( {\frac{5}
{n}} \right)}} {{\Gamma \left( {\frac{1} {n}} \right)}}}} {{\left(
{\sqrt[n]{n}\sigma } \right)^4 \frac{{\Gamma ^2 \left( {\frac{3}
{n}} \right)}} {{\Gamma ^2 \left( {\frac{1} {n}} \right)}}}} =
\frac{{\Gamma \left( {\frac{5} {n}} \right)\Gamma \left( {\frac{1}
{n}} \right)}} {{\Gamma \left( {\frac{3} {n}} \right)\Gamma \left(
{\frac{3} {n}} \right)}}  \hfill \;\;\;\;\;\;\;  (39)  \\
= \frac{{\Gamma \left( {\frac{5} {n}} \right)\Gamma \left( {\frac{1}
{n}} \right)}} {{\Gamma \left( {\frac{3} {n}} \right)\Gamma \left(
{\frac{3} {n}} \right)}}.\frac{{\Gamma \left( {\frac{3} {n} +
\frac{3} {n}} \right)}} {{\Gamma \left( {\frac{5} {n} + \frac{1}
{n}} \right)}} = \frac{{\beta \left( {\frac{1} {n},\frac{5} {n}}
\right)}} {{\beta \left( {\frac{3} {n},\frac{3}
{n}} \right)}} \;\;\;\;\;\;\hfill (40) \\
\end{gathered}\]

And kurtosis excess= $ \frac{{\Gamma \left( {\frac{5} {n}}
\right)\Gamma \left( {\frac{1} {n}} \right)}} {{\Gamma \left(
{\frac{3} {n}} \right)\Gamma \left( {\frac{3} {n}} \right)}} - 3
\hfill (40A) $

According to this the distribution is always mesokurtic, except for
`n=2'.

To be noted that for Normal distribution (n=2), both $m_4$ and $m_2$ are
reducible in terms of `first fundamental moment' i.e. `$ m_0$'(=1).
But for `n=4', `$m_0$'and for `n $>$ 4' neither is expressible in
term of �first fundamental moment�. Hitherto it is assumed that
Normal distribution is the basis of measurement of Skewness,
kurtosis or kurtosis excess of all distributions of practical use.
But keeping the fact in view that, each member of the `family of
curves' represented by $ y(x,n) = \frac{1} {{p_n \sqrt[n]{n}\sigma
}}e^{ - \frac{1} {n}\left( {\frac{{x - m}} {\sigma }} \right)^n }$ ,
(`n' is even number and is the parameter of the family) shows bell
shaped distribution as that of the Normal distribution (the member for n=2). In
that respect, to keep the `analogy', it is rational to express
kurtosis, Kurtosis excess and Coefficient of Skewness in terms of two consecutive
even order moments involving `first fundamental moment' i.e.
`$m_0$'(=1) for any given `n'. In that case Kurtosis, Kurtosis
excess or Skewness values of any distribution obtained from any
experiment subjected to fit normal distribution will depend upon the
order of the normal distribution chosen. Then kurtosis of
normal distribution of order `$\frac {n}{2}$' is $ \frac{{m_{2n} }}
{{m_n^2 }}$.\\
By (36D) $ \frac{{m_{2n} }} {{m_n^2 }} = \frac{{(\sigma ^n )^2
\prod\limits_{r = 0}^{2 - 1} {(1 + rn)} }} {{\left\{ {(\sigma ^n
)\prod\limits_{r = 0}^{1 - 1} {(1 + rn)} } \right\}^2 }} = 1 + n
\hfill  (41)$

Or,

by (38)
\[\begin{gathered} \frac{{m_{2n} }} {{m_n^2 }} = \left( {\sqrt[n]{n}\sigma }
\right)^{2n} \frac{{\Gamma \left( {\frac{1} {n} + \frac{{2n}} {n}}
\right)}} {{\Gamma \left( {\frac{1} {n}} \right)}}/\left\{ {\left(
{\sqrt[n]{n}\sigma } \right)^n \frac{{\Gamma \left( {\frac{1} {n} +
\frac{n} {n}} \right)}} {{\Gamma \left( {\frac{1}
{n}} \right)}}} \right\}^2  \hfill \\
   = \frac{{\Gamma \left( {2 + \frac{1}
{n}} \right).\Gamma \left( {\frac{1} {n}} \right)}} {{\Gamma \left(
{1 + \frac{1} {n}} \right).\Gamma \left( {1 + \frac{1} {n}}
\right)}} = \frac{{\left( {1 + \frac{1} {n}} \right).\frac{1} {n}}}
{{\frac{1} {n}.\frac{1}{n}}} = 1 + n \hfill (41A) \\
\end{gathered}\]

In that respect Kurtosis excess of any obtained distribution is
$\frac{{m_{2n} }} {{m_n^2 }} - (1 + n)$,and that of the chosen normal
distribution itself is `0'.

Coefficient of Skewness of any distribution of form $ y(x,n) = \frac{1} {{p_n
\sqrt[n]{n}\sigma }}e^{ - \frac{1} {n}\left( {\frac{{x - m}} {\sigma
}} \right)^n } $ is `zero', as the curves are symmetric about the
mean of the distribution. Although all odd order moments are zero
for $ y(x,n) = \frac{1} {{p_n \sqrt[n]{n}\sigma }}e^{ - \frac{1}
{n}\left( {\frac{{x - m}} {\sigma }} \right)^n } $ it will be
rational to express Coefficient of Skewness for any given value of n
as, $ \sqrt[n]{{\frac{{m_{n + 1}^n }} {{m_n^{n + 1} }}}}$, in order
to keep the 'analogy' and dimensional similarity (of same order).

\subsection{Multivariate distribution:}

Consider a multivariate distribution for variates $ z_1,\; z_2,\;
z_3,\; ...z_r$. with no correlation, where $ z_i  = \frac{{x_i  -
m_i }}{{\sigma _i }} $. the corresponding orders being $n_1,\;
n_2,\; n_3,\;...n_r$ .then combined P.D.F is given by

$Y_r  = \left[ \frac{1}{{p_{n_1 } p_{n_2 } p_{n_3 } ...p_{n_r }\sqrt[{n1}]{{n1}}\sqrt[{n_2 }]{{n_2 }}\sqrt[{_{n_3 } }]{{n_3 }}...\sqrt[{n_r }]{{n_r }}}} \right] e^{ - \left( {\frac{{z_1^{n_1 }
}}{{n_1 }} + \frac{{z_2^{n_{_2 } } }}{{n_2 }} + \frac{{z_1^{n_3 }}}{{n_3 }} + ... + \frac{{z_r^{n_r } }}{{n_r }}} \right)} $ \\
$\;\;\; = \left[ \frac{1}{{\prod\limits_{i = 1}^r {p_{n_i } \sqrt[{n_i }]{{n_i }}} }}\right] e^{ - \sum\limits_{i = 1}^r {\left( {\frac{{z_i^{n_i } }}{{n_i }}} \right)} }  \hfill (42)$\\

For r =2 it becomes a bi-variate distribution  $ Y_2  = \left[\frac{1}{{p_{n_1 } p_{n_2 } \sqrt[{n1}]{{n1}}\sqrt[{n_2 }]{{n_2 }}}} \right] e^{ - \left( {\frac{{z_1^{n_1 } }}{{n_1 }} +
\frac{{z_2^{n_{_2 } } }}{{n_2 }}} \right)} $

$n_1=\; n_2=2$, it becomes bi-variate Normal distribution. i.e., $ Y_2  = \frac{1}{{2\pi }}e^{ - \frac{1}{2}\left( {z_1^2  + z_2^2 } \right)} \hfill (43)$\\

\section{Extension of Sterling's approximation.}

John Wallis in 1656 [4] has given an expression as follows,
\[
\begin{gathered}
  \left[ {\frac{{2.4.6.8.10....\infty }}
{{1.3.5.7.9....\infty }}} \right]^2  = \frac{\pi }
{2} \;\;\;\;\;\;\;\;\; \hfill (44) \\
  or,\mathop {\lim }\limits_{n \to \infty } \left[ {\frac{{(2n)!!}}
{{(2n - 1)!!}}} \right]^2  = \frac{\pi }
{2} \hfill \\
  or,\mathop {\lim }\limits_{n \to \infty } \left[ {\frac{{(2n)!!}}
{{(2n - 1)!!}}} \right] = \sqrt {\frac{\pi }
{2}}  \hfill \\
  or,\mathop {\lim }\limits_{n \to \infty } (2n)!! = \mathop {\lim }\limits_{n \to \infty } \sqrt {\frac{\pi }
{2}} (2n - 1)!! \hfill \\
  or,\mathop {\lim }\limits_{n \to \infty } \left[ {(2n)!!} \right]^2  = \mathop {\lim }\limits_{n \to \infty } \sqrt {\frac{\pi }
{2}} (2n - 1)!!(2n)!! \;\;\;\;\;\;\;\;\;\;\;\; \hfill (44A) \\
\end{gathered}
\]
Again we have,
\[
\begin{gathered}
  (2n)! = 1.2.3.4.5.6.8.....2n \hfill \\
   = \left( {1.3.5.7.9....2n - 1} \right)\left( {2.4.6.8.10....2n} \right) \hfill \\
   = (2n - 1)!!(2n)!! \hfill \\
  or,(2n - 1)!!(2n)!! = (2n)! \hfill (45) \\
  and, \hfill \\
  (2n)!! = 2.4.6.8.10....2n \hfill \\
   = 2^n (1.2.3.4.5....n) \hfill \\
   = 2^n n! \hfill (46)\\
\end{gathered}
\]
By the equations (44A), (45), (46)...
\[
\begin{gathered}
  \mathop {\lim }\limits_{n \to \infty } \left[ {2^n n!} \right]^2  = \mathop {\lim }\limits_{n \to \infty } \sqrt {\frac{\pi }
{2}} (2n)! \hfill \\
  or,\mathop {\lim }\limits_{n \to \infty } 2^{2n} (n!)^2  = \mathop {\lim }\limits_{n \to \infty } \sqrt {\frac{\pi }
{2}} (2n)! \hfill \\
  or,\mathop {\lim }\limits_{n \to \infty } (2n)! = \mathop {\lim }\limits_{n \to \infty } 2^{2n} \left( {\sqrt {\frac{2}
{\pi }} } \right)(n!)^2  \hfill \\
  or,\mathop {\lim }\limits_{n \to \infty } (2n)! = \mathop {\lim }\limits_{n \to \infty } 2^{2n} \left( {\sqrt {\frac{2}
{\pi }} } \right)\left( {\mathop {\lim }\limits_{n \to \infty } n!} \right)^2  \;\;\;\;\;\;\; \hfill (47)\\
\end{gathered}
\]
From Sterling's approximation [5]
\[
\mathop {\lim }\limits_{n \to \infty } n! \sim \sqrt {2\pi n} \left( {\frac{n}
{e}} \right)^n \;\;\;\;\;\;\;\;\;\;\;\;\;\;\;\;\;\;\;\;\; \hfill (48)
\]
By (47) and (48)
\[
\begin{gathered}
  \mathop {\lim }\limits_{n \to \infty } (2n)! \sim 2^{2n} \left( {\sqrt {\frac{2}
{\pi }} } \right)\left[ {\sqrt {2\pi n} \left( {\frac{n}
{e}} \right)^n } \right]^2  \hfill \\
  or,\mathop {\lim }\limits_{n \to \infty } (2n)! \sim 2^{2n} \left( {\sqrt {\frac{2}
{\pi }} } \right)2\pi n\left( {\frac{n}
{e}} \right)^{2n}  \hfill \\
  or,\mathop {\lim }\limits_{n \to \infty } (2n)! \sim \left( {\sqrt {\frac{2}
{\pi }} } \right)2\pi n\left( {\frac{{2n}}
{e}} \right)^{2n}  \hfill \\
  or,\mathop {\lim }\limits_{n \to \infty } (2n)! \sim 2n\sqrt {2\pi } \left( {\frac{{2n}}
{e}} \right)^{2n} \;\;\;\;\;\;\;\;\;\;\;\;\;\;\;\;\; \hfill (49) \\
\end{gathered}
\]
Again in Sterling's approximation (48), instead of `n' if we put `2n', then,
\[
\begin{gathered}
  \mathop {\lim }\limits_{2n \to \infty } (2n)! \sim \sqrt {2\pi .2n} \left( {\frac{{2n}}
{e}} \right)^{2n}  \hfill \\
  or,\mathop {\lim }\limits_{2n \to \infty } (2n)! \sim \sqrt {2n} \sqrt {2\pi } \left( {\frac{{2n}}
{e}} \right)^{2n}  \;\;\;\;\;\;\;\;\;\;\;\;\;\; \hfill  (50)\\
\end{gathered}
\]
(49)and (50) are two different expressions. Thus Sterling�s approximation can be extended to (49) in order to approximate the value of $ \mathop {\lim }\limits_{n \to \infty } (2n)! $.

\end{document}